\documentclass{article}
\usepackage{graphicx} % Required for inserting images

\usepackage{amsmath,amsfonts,amssymb,amsthm,thmtools,thm-restate,hyperref,color,caption,subcaption,xspace,authblk,fullpage}
\usepackage[shortlabels]{enumitem}
\usepackage[nameinlink]{cleveref}

\usepackage{tikz}
\usetikzlibrary{patterns}
\usetikzlibrary{decorations.pathreplacing}

\newtheorem{theorem}{Theorem}[section]
\newtheorem{lemma}[theorem]{Lemma}

\newtheorem{proposition}[theorem]{Proposition}
\newtheorem{observation}[theorem]{Observation}

\providecommand{\duality}{the \hyperref[lem:duality]{duality lemma}\xspace}
\providecommand{\Duality}{The \hyperref[lem:duality]{duality lemma}\xspace}
\providecommand{\reachability}{the \hyperref[lem:reachability]{reachability lemma}\xspace}
\providecommand{\Reachability}{The \hyperref[lem:reachability]{reachability lemma}\xspace}

\theoremstyle{definition}
\newtheorem{definition}{Definition}[section]

\theoremstyle{remark}

\numberwithin{equation}{section}

\DeclareMathOperator{\row}{\mathsf{row}}
\DeclareMathOperator{\rowone}{\mathsf{row^{=1}}}
\DeclareMathOperator{\rowtwo}{\mathsf{row^{\ge2}}}
\DeclareMathOperator{\col}{\mathsf{col}}
\DeclareMathOperator{\colone}{\mathsf{col^{=1}}}
\DeclareMathOperator{\coltwo}{\mathsf{col^{\ge2}}}
\DeclareMathOperator{\rows}{\mathsf{rows}}
\DeclareMathOperator{\rowsone}{\mathsf{rows^{=1}}}
\DeclareMathOperator{\rowstwo}{\mathsf{rows^{\ge2}}}
\DeclareMathOperator{\cols}{\mathsf{cols}}
\DeclareMathOperator{\colsone}{\mathsf{cols^{=1}}}
\DeclareMathOperator{\colstwo}{\mathsf{cols^{\ge2}}}
\providecommand{\Ztwo}{Z_\ell^{\ge2}}
\providecommand{\Kone}{K^{=1}}
\providecommand{\Ktwo}{K^{\ge2}}
\DeclareMathOperator{\supp}{\mathsf{supp}}
\DeclareMathOperator{\imsupp}{\mathsf{imsupp}}
\providecommand{\RR}{\mathbb{R}}
\DeclareMathOperator{\EE}{\mathbb{E}}

\providecommand{\pmone}{\{-1,1\}}

\providecommand{\tikzZ}{%
\pgfmathsetmacro{\hsep}{0.2}
\pgfmathsetmacro{\wsep}{0.2}
\pgfmathsetmacro{\hh}{1.5}
\pgfmathsetmacro{\ww}{2}
\draw (0,0) -- (3,0) -- (3,-3) -- (0,-3) -- cycle;

\foreach \ii in {1,2,3}
{
 \pgfmathsetmacro{\i}{\ii-1}
 \pgfmathsetmacro{\x}{3+\ii*\wsep+\i*\ww}
 \pgfmathsetmacro{\xx}{3+\ii*\wsep+\ii*\ww}
 \pgfmathsetmacro{\y}{3+\ii*\hsep+\i*\hh}
 \pgfmathsetmacro{\yy}{3+\ii*\hsep+\ii*\hh}
 \draw (\x,-\i) -- (\xx,-\i) -- (\xx,-\ii) -- (\x,-\ii) -- cycle;
 \draw (\i,-\y) -- (\i,-\yy) -- (\ii,-\yy) -- (\ii,-\y) -- cycle;
}

\pgfmathsetmacro{\xleft}{3+\wsep}
\pgfmathsetmacro{\xright}{3+3*(\wsep+\ww)}
\pgfmathsetmacro{\ytop}{3+\hsep}
\pgfmathsetmacro{\ybottom}{3+3*(\hsep+\hh)}
}

\title{Generalized polymorphisms}

\author{Gilad Chase\thanks{E-mail: \texttt{giladchase@gmail.com}}}
\affil{}

\author{Yuval Filmus\thanks{E-mail: \href{mailto:yuvalfi@cs.technion.ac.il}{yuvalfi@cs.technion.ac.il}}}
\affil{The Henry \& Marilyn Taub Faculty of Computer Science \\ Faculty of Mathematics \\ Technion --- Israel Institute of Technology \\ Haifa, Israel}

\begin{document}

\maketitle

\begin{abstract}
    We find all functions $f_0,f_1,\dots,f_m\colon \{0,1\}^n \to \{0,1\}$ and $g_0,g_1,\dots,g_n\colon \{0,1\}^m \to \{0,1\}$ satisfying the following identity for all $n \times m$ matrices $(z_{ij}) \in \{0,1\}^{n \times m}$:
\[
f_0(g_1(z_{11},\dots,z_{1m}),\dots,g_n(z_{n1},\dots,z_{nm})) =
 g_0(f_1(z_{11},\dots,z_{n1}),\dots,f_m(z_{1m},\dots,z_{nm})).
\]
Our results generalize work of Dokow and Holzman (2010), which considered the case $g_0 = g_1 = \cdots = g_n$, and of Chase, Filmus, Minzer, Mossel and Saurabh (2022), which considered the case $g_0 \neq g_1 = \cdots = g_n$.
\end{abstract}

\section{Introduction}
\label{sec:introduction}

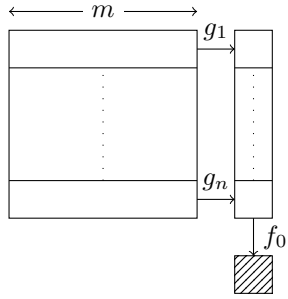
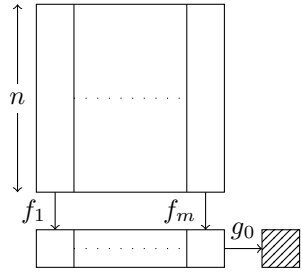
\begin{figure}
\hspace{0.1\textwidth}
\begin{subfigure}{0.35\textwidth}
\centering
\begin{tikzpicture}[scale=0.5]
\draw (0,0) -- (5,0) -- (5,5) -- (0,5) -- cycle;
\draw (0,1) -- (5,1);
\draw (0,4) -- (5,4);
\draw[loosely dotted] (2.5,1) -- (2.5,4);

\draw[<->] (0,5.5) -- (5,5.5) node [midway,fill=white] {$m$};

\draw[->] (5,.5) -- (6,.5) node [midway,above,align=center] {$g_n$};
\draw[->] (5,4.5) -- (6,4.5) node [midway,above,align=center] {$g_1$};

\draw (6,0) -- (7,0) -- (7,5) -- (6,5) -- cycle;
\draw (6,1) -- (7,1);
\draw (6,4) -- (7,4);
\draw[loosely dotted] (6.5,1) -- (6.5,4);

\draw[->] (6.5,0) -- (6.5,-1) node [midway,right,align=center] {$f_0$};

\filldraw[pattern=north east lines] (6,-1) -- (7,-1) -- (7,-2) -- (6,-2) -- cycle;
\end{tikzpicture}
\caption{Left-hand side of \Cref{eq:generalized-polymorphism}}
\end{subfigure}
\hfill
\begin{subfigure}{0.35\textwidth}
\centering
\begin{tikzpicture}[scale=0.5]
\draw (0,0) -- (5,0) -- (5,5) -- (0,5) -- cycle;
\draw (1,0) -- (1,5);
\draw (4,0) -- (4,5);
\draw[loosely dotted] (1,2.5) -- (4,2.5);

\draw[<->] (-0.5,0) -- (-0.5,5) node [midway,fill=white] {$n$};

\draw[->] (0.5,0) -- (0.5,-1) node [midway,left,align=center] {$f_1$};
\draw[->] (4.5,0) -- (4.5,-1) node [midway,left,align=center] {$f_m$};

\draw (0,-1) -- (0,-2) -- (5,-2) -- (5,-1) -- cycle;
\draw (1,-1) -- (1,-2);
\draw (4,-1) -- (4,-2);
\draw[loosely dotted] (1,-1.5) -- (4,-1.5);

\draw[->] (5,-1.5) -- (6,-1.5) node [midway,above,align=center] {$g_0$};

\filldraw[pattern=north east lines] (6,-1) -- (7,-1) -- (7,-2) -- (6,-2) -- cycle;
\end{tikzpicture}
\caption{Right-hand side of \Cref{eq:generalized-polymorphism}}
\end{subfigure}
\hspace{0.1\textwidth}
\caption{Both ways of computing the shaded square should result in the same value}
\label{fig:generalized-polymorphism}
\end{figure}

Boolean functions $f_0,f_1,\ldots,f_m\colon \{0,1\}^n \to \{0,1\}$ and $g_0,g_1,\ldots,g_n\colon \{0,1\}^m \to \{0,1\}$ form a \emph{generalized polymorphism} if the following equation holds for every $0,1$ assignment to the variables $z_{ij}$:
\begin{equation} \label{eq:generalized-polymorphism} f_0(g_1(z_{11},\dots,z_{1m}),\dots,g_n(z_{n1},\dots,z_{nm})) =
 g_0(f_1(z_{11},\dots,z_{n1}),\dots,f_m(z_{1m},\dots,z_{nm})).
\end{equation}
See \Cref{fig:generalized-polymorphism} for an illustration.
In this paper, we give a complete description of all generalized polymorphisms.

We make the following simplifying assumptions, which we later show how to discharge: $f_0,g_0$ depend on all coordinates, and the functions $f_1,\dots,f_m,g_1,\dots,g_n$ are non-constant.

\begin{figure}
\hspace{0.05\textwidth}
\begin{subfigure}{0.4\textwidth}
\centering
\begin{tikzpicture}[scale=0.5]
\tikzZ

\node[align=center] at (1.5,-1.5) {$\Ztwo$};

\draw[decoration={brace},decorate] (0,0.2) -- (3,0.2) node [midway,above] {$\colstwo(Z_\ell)$};
\draw[decoration={brace},decorate] (\xleft,0.2) -- (\xright,0.2) node [midway,above] {$\colsone(Z_\ell)$};
\draw[decoration={brace},decorate] (-0.2,-3) -- (-0.2,0) node [midway,left] {$\rowstwo(Z_\ell)$};
\draw[decoration={brace},decorate] (-0.2,-\ybottom) -- (-0.2,-\ytop) node [midway,left] {$\rowsone(Z_\ell)$};
\end{tikzpicture}
\end{subfigure}
\hfill
\begin{subfigure}{0.4\textwidth}
\centering
\begin{tikzpicture}[scale=0.5]
\tikzZ

\draw[->] (1.3,-1.5) node[below right] {\small{$(i,j)$}} -- (0.5,-0.5);
\draw[dashed] (0,-1) -- (2,-1) -- (2,0);
\draw[dashed] (1,0) -- (1,-2) -- (0,-2);

\draw[decoration={brace},decorate] (0,0.2) -- (2,0.2) node [midway,above] {$\rowtwo(i)$};
\draw[decoration={brace},decorate] (\xleft,0.2) -- (\xright,0.2) node [midway,above] {$\rowone(i)$};
\draw[decoration={brace},decorate] (-0.2,-2) -- (-0.2,0) node [midway,left] {$\coltwo(j)$};
\draw[decoration={brace},decorate] (-0.2,-\ybottom) -- (-0.2,-\ytop) node [midway,left] {$\colone(j)$};
\end{tikzpicture}
\end{subfigure}
\hspace{0.05\textwidth}

\caption{Schematic illustration of the basic definitions for a block $Z_\ell$}
\label{fig:basic}    
\end{figure}
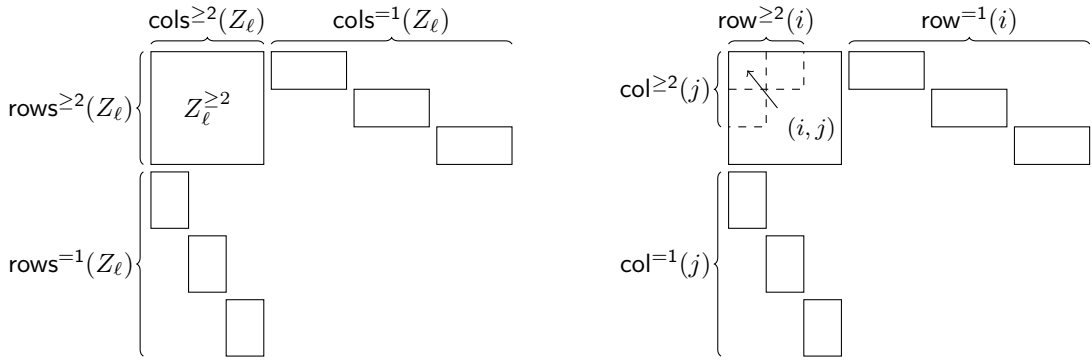

The description of all generalized polymorphisms will require the following terminology, which is illustrated in \Cref{fig:basic}.

\begin{definition} \label{def:basic}
Let $Z \subseteq [n] \times [m]$ be the set of variables that either side of \Cref{eq:generalized-polymorphism} depends on (identifying $z_{ij}$ with $(i,j)$). We think of them as arranged in an $n \times m$ matrix.

Let $\row(i) = \{ j : (i,j) \in Z \}$ be the set of indices on row $i$, which is also the set of coordinates that $g_i$ depends on. Let $\rowone(i) \subseteq \row(i)$ consist of those indices which are unique in their column, and let $\rowtwo(i) = \row(i) \setminus \rowone(i)$. Define $\col(j),\colone(j),\coltwo(j)$ similarly; $\col(j)$ is the set of coordinates that $f_j$ depends on.

We turn $Z$ into a graph by connecting all points on the same row or column. Let $Z_1,\dots,Z_k$ be the corresponding connected components. Let $\cols(Z_\ell)$ be the columns appearing in $Z_\ell$, and define $\colsone(Z_\ell),\colstwo(Z_\ell)$ as above. Similarly, define $\rows(Z_\ell),\rowsone(Z_\ell),\rowstwo(Z_\ell)$. Finally, let $\Ztwo$ consist of those $(i,j) \in Z_\ell$ such that $i \in \rowstwo(Z_\ell)$ and $j \in \colstwo(Z_\ell)$.
\end{definition}

\begin{theorem}[Main] \label{thm:main}
Let $f_0,f_1,\dots,f_m\colon \{0,1\}^n \to \{0,1\}$ and $g_0,g_1,\dots,g_n\colon \{0,1\}^m \to \{0,1\}$ be a generalized polymorphism, where $f_0,g_0$ depend on all coordinates, and the functions $f_1,\dots,f_m,g_1,\dots,g_n$ are non-constant.

Let $Z$ be the set of variables that either side of \Cref{eq:generalized-polymorphism} depends on, and let $Z_1,\dots,Z_k$ be the connected components of the corresponding graph.

There exist a function $h\colon \{0,1\}^k \to \{0,1\}$ and, for each $\ell \in [k]$, functions $f_0^{(\ell)}\colon \{0,1\}^{\rows(Z_\ell)} \to \{0,1\}$ and $g_0^{(\ell)}\colon \{0,1\}^{\cols(Z_\ell)} \to \{0,1\}$, such that
\begin{align*}
 f_0(x) &= h\bigl(f_0^{(1)}(x|_{\rows(Z_1)}), \dots, f_0^{(k)}(x|_{\rows(Z_k)})\bigr), &
 g_0(y) &= h\bigl(g_0^{(1)}(y|_{\cols(Z_1)}), \dots, g_0^{(k)}(y|_{\cols(Z_k)})\bigr).
\end{align*}

For each $\ell \in [k]$, the functions $f_0^{(\ell)},g_0^{(\ell)}$, the functions $f_j$ for $j \in \cols(Z_\ell)$, and the functions $g_i$ for $i \in \rows(Z_\ell)$, are as follows.

If $j \in \colsone(Z_\ell)$ and $(i,j) \in Z_\ell$ then for some $\sigma_j \in \{0,1\}$, $f_j(x) = x_i \oplus \sigma_j$.

If $i \in \rowsone(Z_\ell)$ and $(i,j) \in Z_\ell$ then for some $\tau_i \in \{0,1\}$, $g_i(y) = y_j \oplus \tau_i$.

For the remaining functions, we have three different possibilities.

\paragraph{Singleton case} $Z_\ell = \{(i,j)\}$, $f_0^{(\ell)}(x) = x_i \oplus \tau_i$, and $g_0^{(\ell)}(y) = y_j \oplus \sigma_j$.

\paragraph{XOR case} There are functions $\gamma_i\colon \{0,1\}^{\rowone(i)} \to \{0,1\}$ for each $i \in \rowstwo(Z_\ell)$ and $\phi_j\colon \{0,1\}^{\colone(j)} \to \{0,1\}$ for each $j \in \colstwo(Z_\ell)$ such that
\begin{align*}
 f_0^{(\ell)}(x) &= \bigoplus_{i \in \rowstwo(Z_\ell)} x_i \oplus \bigoplus_{j \in \colstwo(Z_\ell)} \phi_j((x \oplus \tau)|_{\colone(j)}), &
 g_i(y) &= \bigoplus_{j \in \rowtwo(i)} y_j \oplus \gamma_i(y|_{\rowone(i)}), \\
 g_0^{(\ell)}(y) &= \bigoplus_{j \in \colstwo(Z_\ell)} y_j \oplus \bigoplus_{i \in \rowstwo(Z_\ell)} \gamma_i((y \oplus \sigma)|_{\rowone(i)}), &
 f_j(x) &= \bigoplus_{i \in \coltwo(j)} x_i \oplus \phi_j(x|_{\colone(j)}).
\end{align*}

\paragraph{AND case} There are functions $\gamma_i,\phi_j$ as in the XOR case;
%There are functions $\gamma_i\colon \{0,1\}^{\rowone(i)} \to \{0,1\}$ for each $i \in \rowstwo(Z_\ell)$ and $\phi_j\colon \{0,1\}^{\colone(j)} \to \{0,1\}$ for each $j \in \colstwo(Z_\ell)$;
constants $D_i \in \{0,1\}$ for each $i \in \rowstwo(Z_\ell)$ and $B_j \in \{0,1\}$ for each $j \in \colstwo(Z_\ell)$; and constants $\kappa_{ij} \in \{0,1\}$ for all $(i,j) \in \Ztwo$, such that
\begin{align*}
 f_0^{(\ell)}(x) &= \bigwedge_{i \in \rowstwo(Z_\ell)} (x_i \oplus D_i) \land \bigwedge_{j \in \colstwo(Z_\ell)} \phi_j((x \oplus \tau)|_{\colone(j)}), &
 g_i(y) &= \left(\bigwedge_{j \in \rowtwo(i)} (y_j \oplus \kappa_{ij}) \land \gamma_i(y|_{\rowone(i)})\right) \oplus D_i, \\
 g_0^{(\ell)}(y) &= \bigwedge_{j \in \colstwo(Z_\ell)} (y_j \oplus B_j) \land \bigwedge_{i \in \rowstwo(Z_\ell)} \gamma_i((y \oplus \sigma)|_{\rowone(i)}), &
 f_j(x) &= \left(\bigwedge_{i \in \coltwo(j)} (x_i \oplus \kappa_{ij}) \land \phi_j(x|_{\colone(j)})\right) \oplus B_j.
\end{align*}
\end{theorem}

Conversely, every collection of functions $f_0,f_1,\dots,f_m\colon \{0,1\}^n \to \{0,1\}$ and $g_0,g_1,\dots,g_n\colon \{0,1\}^m \to \{0,1\}$ which have the form given in the main theorem constitute a generalized polymorphism. To see this, it suffices to verify that for all $\ell \in [k]$,
\[
 f_0^{(\ell)}(g_i(z_{i1},\dots,z_{im}) : i \in \rows(Z_\ell)) =
 g_0^{(\ell)}(f_j(z_{1j},\dots,z_{nj}) : j \in \cols(Z_\ell)).
\]
In the singleton case, the common value is $z_{ij}$, and in the other two cases, it is one of
\begin{gather*}
\bigoplus_{(i,j) \in \Ztwo} z_{ij} \oplus
\bigoplus_{i \in \rowstwo(Z_\ell)} \gamma_i(z_{ij} : j \in \rowone(i)) \oplus
\bigoplus_{j \in \colstwo(Z_\ell)} \phi_j(z_{ij} : i \in \colone(j)), \\
\bigwedge_{(i,j) \in \Ztwo} z_{ij} \land
\bigwedge_{i \in \rowstwo(Z_\ell)} \gamma_i(z_{ij} : j \in \rowone(i)) \land
\bigwedge_{j \in \colstwo(Z_\ell)} \phi_j(z_{ij} : i \in \colone(j)).
\end{gather*}

The following observation allows us to discharge the simplifying assumptions.

\begin{observation} \label{obs:simplifying-assumptions}
Let $f_0,f_1,\dots,f_m\colon \{0,1\}^n \to \{0,1\}$ and $g_0,g_1,\dots,g_n\colon \{0,1\}^m \to \{0,1\}$ be an arbitrary generalized polymorphism.

Let $I$ be the set of $i \in [n]$ such that $g_i$ is non-constant, and let $J$ be the set of $j \in [m]$ such that $f_j$ is non-constant.

Let $f_0|_I\colon \{0,1\}^I \to \{0,1\}$ be obtained from $f_0$ by substituting the constant value of $g_i$ for the coordinates $i \notin I$, and let $I' \subseteq I$ be the set of coordinates that $f_0$ depends on. Define $g_0|_J$ and $J'$ similarly.

If $i \in I'$ then $g_i$ is a non-constant function which only depends on coordinates in $J'$.
\end{observation}
\begin{proof}
If $i \in I'$ and $g_i$ depends on $j$ then \Cref{eq:generalized-polymorphism} depends on $z_{ij}$. This can only happen if $j \in J'$. Since $i \in I$, the function $g_i$ is non-constant.
\end{proof}

By definition, $f_0|_I$ depends only on the coordinates in $I'$. Let $f_0|_{I'}\colon \{0,1\}^{I'} \to \{0,1\}$ be the corresponding restriction of $f_0|_I$. The observation shows that for $i \in I'$, the function $g_i$ depends only on the coordinates in $J'$. Let $g_i|_{J'}\colon \{0,1\}^{J'} \to \{0,1\}$ be the corresponding restriction of $g_i$. Restricting the duals of these functions, we obtain a generalized polymorphism
\[
 f_0|_{I'}, (f_j|_{I'})_{j \in J'}, g_0|_{J'}, (g_i|_{J'})_{i \in I'}
\]
in which $f_0|_{I'},g_0|_{J'}$ depend on all coordinates and $f_j|_{I'},g_i|_{J'}$ are non-constant. \Cref{thm:main} describes the structure of this restricted generalized polymorphism, and hence of the original one.

\subsection{Background} \label{sec:background}

\begin{figure}
\centering
\begin{tikzpicture}[scale=0.5]
\draw (0,0) -- (5,0) -- (5,5) -- (0,5) -- cycle;
\draw (1,0) -- (1,5);
\draw (4,0) -- (4,5);
\draw[loosely dotted] (1,2.5) -- (4,2.5);

\draw[<->] (0,5.5) -- (5,5.5) node [midway,fill=white] {$p$};
\draw[<->] (-0.5,0) -- (-0.5,5) node [midway,fill=white] {$n$};

\draw[->] (0.5,0) -- (0.5,-1) node [midway,left,align=center] {$f$};
\draw[->] (4.5,0) -- (4.5,-1) node [midway,left,align=center] {$f$};

\draw (0,-1) -- (0,-2) -- (5,-2) -- (5,-1) -- cycle;
\draw (1,-1) -- (1,-2);
\draw (4,-1) -- (4,-2);
\draw[loosely dotted] (1,-1.5) -- (4,-1.5);

\node [right] at (5,0.5) {$\in P$};
\draw[loosely dotted] (5.75,1) -- (5.75,4);
\node [right] at (5,4.5) {$\in P$};

\draw[double distance=1.2pt,->] (5.75,0) -- (5.75,-1);

\node [right] at (5,-1.5) {$\in P$};

\end{tikzpicture}

\caption{If all rows of the $n \times p$ table belong to $P$, then so does the additional row}
\label{fig:polymorphism}
\end{figure}
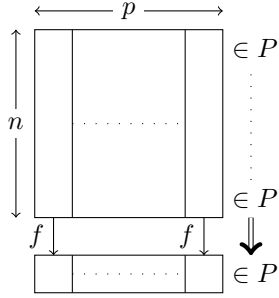

Polymorphisms make an appearance in two different areas of mathematics: universal algebra (where the name originates) and social choice theory. Typically, one considers polymorphisms of \emph{predicates}. Given a finite alphabet $\Sigma$ and a $p$-ary predicate $P \subseteq \Sigma^p$, a function $f\colon \Sigma^n \to \Sigma$ is a polymorphism of $P$ if given $(z_{i1},\dots,z_{ip}) \in P$ for $i \in [n]$, also $(f(z_{11},\dots,z_{n1}), \dots, f(z_{1p},\dots,f_{np})) \in P$. See \Cref{fig:polymorphism} for an illustration. Polymorphisms of predicates are used in the modern approach to Schaefer's dichotomy theorem~\cite{Chen09} and its generalizations~\cite{Bulatov2017,Zhuk2020}.

In social choice theory, the predicate $P$ is called a \emph{domain} or an \emph{agenda}. A domain describes the set of consistent \emph{judgments}. For example, if voters rank their preference over three candidates $A,B,C$ by stating which candidate they prefer out of each pair $(A,B),(B,C),(C,A)$, then the result satisfies the predicate $P_{NAE} \subseteq \{0,1\}^3$ consisting of all triplets which are not all equal (since we cannot have $A \prec B \prec C \prec A$ or $A \succ B \succ C \succ A$).
The opinions of the voters can be aggregated using a judgment aggregation function. A judgment aggregation function is \emph{consistent} if it is a polymorphism of the domain.
Arrow's classical result~\cite{Arrow} states that the only polymorphisms of $P_{NAE}$ which satisfy \emph{unanimity} ($f(b,\ldots,b) = b$) are dictators (functions of the form $f(x) = x_i$).

A domain is \emph{truth-functional} if it has the form $P_g = \{(y_1,\ldots,y_{m+1}) : y_{m+1} = g(y_1,\ldots,y_m)\}$. For example, the classical work of Kornhauser and Sager~\cite{KornhauserSager} highlights the fact that Majority is not a polymorphism of the predicate $P_{\land} = \{(a,b,c) \in \{0,1\}^3 : a \land b = c\}$.
A function $f$ is a polymorphism of the predicate $P_g$ if $f_0 = \cdots = f_m = f$ and $g_0 = \cdots = g_n = g$ constitute a generalized polymorphism per \Cref{eq:generalized-polymorphism}.

From the point of view of social choice theory, it is natural to allow the judgment aggregation function to depend on the coordinate; as shown by Szegedy and Xu~\cite{SX15,Xu15}, this corresponds to the multisorted setting in universal algebra~\cite{BulatovJeavons2003,Bulatov2011}. Dokow and Holzman~\cite{DH09}, motivated by social choice theory, consider this more general setting. In our terminology, they describe all generalized polymorphisms (over the Boolean alphabet) in which $g_0 = \cdots = g_n = g$. Chase et al.~\cite{CFMMS22} extend this to the setting in which $g_1 = \cdots = g_n = g$, but possibly $g_0 \neq g$.

Dokow and Holzman~\cite{DH10} consider generalized polymorphisms (under the assumption $g_0 = \cdots = g_n = g$) for arbitrary finite alphabets. They determine for which functions $g$, the only polymorphisms are trivial (constant functions, dictators, or negations of dictators).

A completely different line of work deals with approximate polymorphisms, which are functions $f$ satisfying the polymorphism condition for most input matrices. Most of the work in this direction has dealt with Arrow's theorem and its generalizations~\cite{Kalai,Keller,Mossel2012,FKKN,IKM,MR,FalikFriedgut}. Nehama~\cite{Nehama} and Filmus et al.~\cite{FLMM2020} consider $P_{\land}$, and Chase et al.~\cite{CFMMS22} consider arbitrary truth-functional predicates over $\{0,1\}$.

\subsection{Paper organization} \label{sec:paper-organization}

The rest of the paper constitutes the proof of \Cref{thm:main}. To avoid trivialities, we assume throughout that $n,m \ge 1$, an assumption we don't explicitly mention. In order to avoid clutter, we assume an additional simplifying assumption: if any of the functions $f_1,\dots,f_m,g_1,\dots,g_n$ depends on a single coordinate, then it is a dictator (of the form $f(x) = x_i$) rather than an anti-dictator (of the form $f(x) = \lnot x_i$).

The proof follows the approach used in Chase et al.~\cite{CFMMS22}, which uses Fourier analysis. For this reason, we switch from $\{0,1\}$ to $\{-1,1\}$. That is, a generalized polymorphism now consists of functions $f_0,f_1,\dots,f_m\colon \pmone^n \to \pmone$ and functions $g_0,g_1,\dots,g_n\colon \pmone^m \to \pmone$ which satisfy \Cref{eq:generalized-polymorphism}.

Let us briefly recall some standard notions in Fourier analysis. Each function $f\colon \pmone^n \to \pmone$ has a unique expansion as a multilinear polynomial, known as its \emph{Fourier expanson}:
\[
 f(x) = \sum_{S \subseteq [n]} \hat{f}(S) \prod_{i \in S} x_i.
\]
The coefficients $\hat{f}(S)$ are known as \emph{Fourier coefficients}. The \emph{Fourier support} of $f$ is $\supp(\hat{f}) = \{ S \subseteq [n] : \hat{f}(S) \neq 0 \}$. %The \emph{degree} of $f$, denoted $\deg f$, is the maximum size of a set in $\supp(\hat{f})$, or zero if $f$ is constant.

A set $S$ in the Fourier support of $f$ is \emph{inclusion-maximal} if no set in the Fourier support properly contains $S$. We denote the family of inclusion-maximal sets in the Fourier support by $\imsupp(\hat{f})$.\footnote{This notation is non-standard. We are not aware of any standard notation for the inclusion-maximal support.}

% try to get away without using degree? otherwise, define it here

The major step in the proof is \emph{shifting}: we find parameters $\kappa_{ij},B_j,D_i$ such that the functions
\begin{align*}
    G_i(y) &= g_i(y_1 - \kappa_{i1}, \dots, y_m - \kappa_{im}), &
    F_j(x) &= f_j(x_1 - \kappa_{1j}, \dots, x_n - \kappa_{jn})
\end{align*}
have the particularly simple form
\begin{align*}
    G_i(y) &= \prod_{j \in \rowtwo(i)} y_j \cdot \gamma_i(y|_{\rowone(i)}), &
    F_j(x) &= \prod_{i \in \coltwo(j)} x_i \cdot \phi_j(x|_{\colone(j)}).
\end{align*}
In particular, they are \emph{centered}: they have zero mean. %Moreover, they are all \emph{Boolean-like}, a concept we define in \Cref{sec:preliminaries}.

If we define
\begin{align*}
    F_0(x) &= f_0(x_1 + D_1, \dots, x_n + D_n), &
    G_0(y) &= g_0(y_1 + B_1, \dots, y_m + B_m),
\end{align*}
then the functions $F_0,F_1,\dots,F_m,G_0,G_1,\dots,G_n$ satisfy \Cref{eq:generalized-polymorphism}, which we use in the proof. Note that these functions are not necessarily Boolean, that is, they might not map $\pmone^n$ or $\pmone^m$ to $\pmone$.

Here is how the remainder of the paper is organized:
\begin{itemize}
\item \Cref{sec:preliminaries} contains some basic definitions, and two important lemmas: the duality lemma and the reachability lemma. The duality lemma is a property of the Fourier support of centered functions, and the reachability lemma is a reachability property for the connected components $Z_\ell$.
\item \Cref{sec:step1,sec:step2} prove properties of the inclusion-maximal Fourier support of the functions $F_0,G_0$ (\Cref{sec:step1}) and $F_1,\dots,F_n,G_1,\dots,G_m$ (\Cref{sec:step2}), only assuming that they are non-constant, centered, and Boolean-like (a property defined in \Cref{sec:preliminaries} which holds for all shifts of Boolean functions).
\item \Cref{sec:step3} performs the shifting step, and proves that $F_1,\dots,F_m,G_1,\dots,G_n$ have the stated structure. \Cref{sec:step4} then goes on to describe the structure of $F_0,G_0$.
\item \Cref{sec:step5} deduces \Cref{thm:main} by first undoing the shifting, and then using a structure theorem for functions of the form
\[
 f(u,v) = \phi(u) \prod_i (v_i + \kappa_i) + B,
\]
which we prove in \Cref{sec:step5-prel}.
\end{itemize}

The definition of generalized polymorphism is symmetric in rows and columns. Accordingly, many statements arising during the proof have duals (for example, a statement about $\imsupp(\hat{f}_0)$ has a dual about $\imsupp(\hat{g}_0)$). In such cases, we only prove one of these statements, but freely use both subsequently.

\paragraph{Acknowledgements}
This paper constitutes the first author's MSc thesis.
We thank the thesis committee for their helpful remarks.

This project has received funding from the European Union's Horizon 2020 research and innovation programme under grant agreement No~802020-ERC-HARMONIC.

\section{Preliminaries} \label{sec:preliminaries}

Throughout the paper, we identify functions over $\pmone^n$ with multilinear polynomials over $n$ variables (this is possible since every function over $\pmone^n$ has a unique representation as a multilinear polynomial). A function $f$ over $\pmone^n$ is \emph{Boolean} if $f(x) \in \pmone$ for every $x \in \pmone^n$. A function $f$ is \emph{centered} if $\hat{f}(\emptyset) = 0$.

In the introduction, we defined a generalized polymorphism to be a set of Boolean functions satisfying \Cref{eq:generalized-polymorphism}. As explained in \Cref{sec:paper-organization}, part of the proof will involve a set of non-Boolean functions satisfying \Cref{eq:generalized-polymorphism}. For succinctness, we will call the latter a generalized polymorphism, and the former a Boolean generalized polymorphism.

\begin{definition}[Generalized polymorphism]
Let $f_0,f_1,\dots,f_m$ be $n$-ary multilinear polynomials, and let $g_0,g_1,\dots,g_n$ be $m$-ary multilinear polynomials. The functions $f_0,f_1,\dots,f_m,g_0,g_1,\dots,g_n$ constitute a \emph{generalized polymorphism} if they satisfy \Cref{eq:generalized-polymorphism}. If all of these functions are also Boolean, then they constitute a \emph{Boolean generalized polymorphism}.
\end{definition}

As explained in \Cref{sec:paper-organization}, non-Boolean generalized polymorphisms arise by shifting the original generalized polymorphism.

\begin{definition}[Shift] \label{def:shift}
Let $f\colon \pmone^n \to \RR$. A \emph{shift} of $f$ is any function $F\colon \pmone^n \to \RR$ given by
\[
 F(x) = f(x_1 - \kappa_1, \dots, x_n - \kappa_n) - B,
\]
where $\kappa_1,\dots,\kappa_n,B \in \mathbb{R}$.
\end{definition}

Shifting preserves many crucial properties.

\begin{observation} \label{obs:shifting}
If $F$ is a shift of a non-constant function $f$ then $f$ and $F$ depend on the same coordinates, and moreover $\imsupp(f) = \imsupp(F)$.
\end{observation}

As a consequence, the connected components $Z_1,\dots,Z_k$ and the various other definitions appearing in \Cref{sec:introduction} are the same for the original generalized polymorphism $f_0,f_1,\dots,f_m,g_0,g_1,\dots,g_n$ and for the shifted generalized polymorphism $F_0,F_1,\dots,F_m,G_0,G_1,\dots,G_n$ mentioned in \Cref{sec:paper-organization}.

\paragraph{Boolean-like functions}

\Cref{sec:step1,sec:step2} will need to assume that all shifted polymorphisms satisfy the following property.

% let's hope that this definition works...

\begin{definition}[Boolean-like] \label{def:boolean-like}
A function $F\colon \pmone^n \to \pmone$ is \emph{Boolean-like} if the following holds: either (i) $F$ depends on at most one coordinate, or (ii) all non-empty $S \in \imsupp(\hat{F})$ satisfy $|S| \ge 2$.
\end{definition}

\begin{observation} \label{obs:boolean-like}
Every Boolean function is Boolean-like, and every shift of a non-constant Boolean-like function is a non-constant Boolean-like function.
\end{observation}
\begin{proof}
The second property follows immediately from the observation that shifting preserves $\imsupp$, so we only prove the first property.

Let $f\colon \pmone^n \to \pmone$, and suppose that $\imsupp(\hat{f})$ contains a singleton, say $\{x_1\}$. Then we can write $f(x) = A + Bx_1 + g(x_2,\dots,x_n)$. Assume, without loss of generality, that $B > 0$. Clearly $f(1,x_2,\dots,x_n) > f(-1,x_2,\dots,x_n)$, and so $f(x_1,x_2,\dots,x_n) = x_1$. Therefore $g$ must be constant, and so $f$ depends only on $x_1$.
\end{proof}

\subsection{Duality lemma} \label{sec:duality}

When the functions $f_1,\dots,f_m,g_1,\dots,g_n$ are centered, there is a certain duality between the Fourier supports. Let $W \subseteq [n] \times [m]$. We can describe $W$ in two ways:
\begin{itemize}
    \item The set of row indices appearing in $W$: $\rows(W) = \{ i : (i,j) \in W \text{ for some } j \}$. \\
    For each $i \in \rows(W)$, the corresponding columns: $\row(W,i) = \{ j : (i,j) \in W \}$.
    \item The set of column indices appearing in $W$: $\cols(W) = \{ j : (i,j) \in W \text{ for some } i \}$. \\
    For each $j \in \cols(W)$, the corresponding rows: $\col(W,j) = \{ i : (i,j) \in W \}$.
\end{itemize}

\begin{lemma}[Duality lemma] \label{lem:duality}
Suppose $f_0,f_1,\dots,f_m,g_0,g_1,\dots,g_n$ is a generalized polymorphism, where $f_1,\dots,f_m,g_1,\dots,g_n$ are centered and non-zero.

Let $W \subseteq [n] \times [m]$. Then:
\begin{enumerate}[(a)]
\item $\rows(W) \in \supp(f_0)$ and $\row(W,i) \in \supp(g_i)$ for all $i \in \rows(W)$ if and only if \\
$\cols(W) \in \supp(g_0)$ and $\col(W,j) \in \supp(f_j)$ for all $j \in \cols(W)$.
\item $\rows(W) \in \imsupp(f_0)$ and $\row(W,i) \in \imsupp(g_i)$ for all $i \in \rows(W)$ if and only if \\
$\cols(W) \in \imsupp(g_0)$ and $\col(W,j) \in \imsupp(f_j)$ for all $j \in \cols(W)$.
\end{enumerate}
\end{lemma}
\begin{proof}
Substituting the Fourier expansions of the various functions in \Cref{eq:generalized-polymorphism}, we obtain
\[
 \sum_{\substack{A \subseteq [n] \\ \text{for all } i \in A\colon \\ B_i \subseteq [m], B_i \neq \emptyset}}
 \hat{f}_0(A) \prod_{i \in A} \hat{g}_i(B_i) \prod_{j \in B_i} z_{ij}
 =
 \sum_{\substack{B \subseteq [m] \\ \text{for all } j \in B\colon \\ A_j \subseteq [n], A_j \neq \emptyset}}
 \hat{g}_0(B) \prod_{j \in B} \hat{f}_j(A_j) \prod_{i \in A_j} z_{ij},
\]
where we can assume that $B_i,A_j \neq \emptyset$ since $f_1,\dots,f_m,g_1,\dots,g_n$ are centered.
Since $B_i,A_j \neq \emptyset$, every monomial appears exactly once on each side of this equation.

Extracting the coefficient of the monomial $\prod_{(i,j) \in W} z_{ij}$, we obtain
\[
 \hat{f}_0(\rows(W)) \prod_{i \in \rows(W)} \hat{g}_i(\row(W,i)) = 
 \hat{g}_0(\cols(W)) \prod_{j \in \cols(W)} \hat{f}_j(\col(W,j)).
\]
This immediately implies part~(a).

\medskip

We prove part~(b) by contradiction. Suppose that $\rows(W) \in \imsupp(f_0)$ and $\row(W,i) \in \imsupp(g_i)$ for all $i \in \rows(W)$. We first show that $\col(W,j)$ is inclusion-maximal for all $j \in \cols(W)$, and then that $\cols(W)$ itself is inclusion-maximal.

Suppose, for the sake of contradiction, that $\col(W,j_0)$ is not inclusion-maximal, for some $j_0 \in \cols(W)$. Then there exists a set $S \supsetneq \col(W,j_0)$ in the support of $\hat{f}_{j_0}$. Let $W'$ be obtained by adding to $W$ the points $(i,j_0)$ for $i \in S \setminus \col(W,j_0)$, and choose an arbitrary $i_0 \in S \setminus \col(W,j_0)$.

We now consider two cases. If $i_0 \notin \rows(W)$ then applying part~(a) (in the $\gets$ direction) to $W'$, we deduce that $\rows(W') = \rows(W) \cup \{i_0\} \in \supp(\hat{f}_0)$, contradicting the inclusion-maximality of $\rows(W)$. If $i_0 \in \rows(W)$ then applying part~(a), we deduce that $\row(W',i_0) = \row(W,i_0) \cup \{j_0\} \in \supp(\hat{g}_{i_0})$, contradicting the inclusion-maximality of $\row(W,i_0)$.

\smallskip

Suppose now that $\cols(W)$ is not inclusion-maximal. Let $T \supsetneq \cols(W)$ be in the support of $\hat{g}_0$. For each $j \in T \setminus \cols(W)$, choose an arbitrary set $S_j \in \imsupp(\hat{f}_j)$ (this is possible since $f_j$ is non-zero). Let $W'$ be obtained by adding to $W$ the points $(i,j)$ for each $j \in T \setminus \cols(W)$ and $i \in S_j$.

Apply part~(a) (in the $\gets$ direction) to $W'$ to deduce that $\rows(W') = \rows(W) \cup \bigcup_{j \in T \setminus \cols(W)} S_j \in \supp(\hat{f}_0)$. By inclusion-maximality, $\rows(W') = \rows(W)$, and so $S_j \subseteq \rows(W)$ for all $j \in T \setminus \cols(W)$.

Part~(a) also implies that for every $i \in \rows(W)$, we have $\row(W',i) = \row(W,i) \cup \{ j \in T \setminus \cols(W) : i \in S_j \} \in \supp(\hat{g}_i)$. By inclusion-maximality, $\row(W',i) = \row(W,i)$, and so $S_j$ must be disjoint from $\rows(W)$ for all $j \in T \setminus \cols(W)$.

Put together, we see that $S_j = \emptyset$ for all $j \in T \setminus \cols(W)$. Since $f_j$ is centered, this is impossible.
\end{proof}

\subsection{Reachability lemma} \label{sec:reachability}

Let $f_0,f_1,\dots,f_m,g_0,g_1,\dots,g_n$ be a generalized polymorphism. If $\emptyset \subsetneq I \subsetneq \rows(Z_\ell)$ then since $Z_\ell$ is a connected component, we can find $i_0 \in I$ and $i_1 \in \rows(Z_\ell) \setminus I$ that appear on the same column $j \in \cols(Z_\ell)$, that is, $(i_0,j),(i_1,j) \in Z_\ell$. We will need a refinement of this observations.

\begin{lemma}[Reachability lemma] \label{lem:reachability}
Let $f_0,f_1,\dots,f_m,g_0,g_1,\dots,g_n$ be a generalized polymorphism, with connected components $Z_1,\dots,Z_k$.

Let $\ell \in [k]$, and let $I \subseteq [n]$ be a set which intersects $\rows(Z_\ell)$ but doesn't contain $\rowstwo(Z_\ell)$. Then there exist $i_0 \in I \cap \rows(Z_\ell)$ and $i_1 \in \rowstwo(Z_\ell) \setminus I$ such that $i_0,i_1$ appear on the same column $j$ for some $j \in \cols(Z_\ell)$, that is $(i_0,j),(i_1,j) \in Z_\ell$.
\end{lemma}
\begin{proof}
We can assume, without loss of generality, that $I \subseteq \rows(Z_\ell)$.
Extend $I$ to the set
\[
 I' = I \cup \{ i' \in \rowsone(Z_\ell) : \text{$i'$ appears on the same column as some $i \in I$} \}.
\]
Since $I$ doesn't contain all of $\rowstwo(Z_\ell)$, the same holds for $I'$, and in particular, $\emptyset \subsetneq I' \subsetneq \rows(Z_\ell)$. Since $Z_\ell$ is a connected component, there must exist $i'_0 \in I'$ and $i'_1 \in \rows(Z_\ell) \setminus I'$ that appear on the same column $j$ for some $j \in \cols(Z_\ell)$.

We consider two cases: $i'_0 \in I$ and $i'_0 \notin I$. Suppose first that $i'_0 \in I$. We claim that $i'_1 \in \rowstwo(Z_\ell)$. Indeed, if $i'_1 \in \rowsone(Z_\ell)$ then $i'_1$ appears on the same column as $i'_0 \in I$, and so $i'_1 \in I'$, contradicting the definition of $i'_1$. We can therefore take $i_0 = i'_0$ and $i_1 = i'_1$.

Suppose next that $i'_0 \notin I$. By definition of $I'$, we have $i'_0 \in \rowsone(Z_\ell)$, and there exists $i''_0 \in I$ such that $i'_0$ and $i''_0$ are on the same column, which must be column $j$ since $\row(i'_0) = \{j\}$. As in the preceding case, $i'_1 \in \rowstwo(Z_\ell)$, and so we can take $i_0 = i''_0$ and $i_1 = i'_1$.
\end{proof}

\section{Inclusion-maximal support of \texorpdfstring{$f_0,g_0$}{outer functions}}
\label{sec:step1}

In this section we prove the following result.

\begin{proposition} \label{pro:step1}
Let $F_0,F_1,\dots,F_m,G_0,G_1,\dots,G_n$ be a generalized polymorphism, where $F_0,G_0$ depend on all coordinates, and $F_1,\dots,F_m,G_1,\dots,G_n$ are Boolean-like, centered, and non-zero.

Let $Z_1,\dots,Z_k$ be the corresponding connected components, and let $\ell \in [k]$.

If $S \in \imsupp(\hat{f}_0)$ intersects $\rows(Z_\ell)$ then $S$ contains $\rowstwo(Z_\ell)$.
\end{proposition}
\begin{proof}
The proof is in two parts. In the first part, we show that $\supp(\hat{f}_0)$ includes a set $S_0$ which contains $\rowstwo(Z_\ell)$ and agrees with $S$ outside of $\rows(Z_\ell)$ (that is, $S \setminus \rows(Z_\ell) = S_0 \setminus \rows(Z_\ell)$). In the second part, we use $S_0$ to show that $S$ itself must contain $\rowstwo(Z_\ell)$.

Throughout the proof we assume that $\rowstwo(Z_\ell)$ is non-empty, since otherwise the \namecref{pro:step1} becomes trivial.

\paragraph{First part} Among all sets in $\supp(\hat{f}_0)$ which intersect $\rows(Z_\ell)$ and agree with $S$ outside of $\rows(Z_\ell)$, let $S_0$ maximize $|S_0 \cap \rowstwo(Z_\ell)|$. We will show that in fact, $S_0 \supseteq \rowstwo(Z_\ell)$.

Assume, for the sake of contradiction, that $S_0$ doesn't contain $\rowstwo(Z_\ell)$. Since $S_0$ intersects $\rows(Z_\ell)$, \reachability shows that there exist $i_0 \in S_0$ and $i_1 \in \rowstwo(Z_\ell) \setminus S_0$ such that $i_0,i_1$ appear on the same column $j_0 \in \cols(Z_\ell)$.

Define $W_0 \subseteq [n] \times [m]$ as follows: $\rows(W_0) = S_0$; $\row(W_0,i_0)$ is some set in $\imsupp(\hat{g}_{i_0})$ containing $j_0$ (this is possible since $(i_0,j_0) \in Z$); and $\row(W,i)$ is an arbitrary set in $\imsupp(\hat{g}_i)$ for $i \in S_0$ other than $i_0$.
Applying \duality, we deduce that $\cols(W_0) \in \supp(\hat{g}_0)$ and $\col(W_0,j) \in \supp(\hat{f}_j)$ for all $j \in \cols(W_0)$.
Note $(i_0,j_0) \in W_0$.

Let $W_1$ be obtained from $W_0$ by replacing column $j_0$ with some set in $\imsupp(\hat{f}_{j_0})$ containing $i_1$; this is possible since $(i_1,j_0) \in Z$. Applying \duality again, we deduce that $\rows(W_1) \in \supp(\hat{f}_0)$.

Since $S_0$ and $\rows(W_1)$ differ only on rows appearing in column $j_0 \in \cols(Z_\ell)$, it follows that $S_0$ and $\rows(W_1)$ agree outside of $\rows(Z_\ell)$.
We will show that moreover, $\rows(W_1) \cap \rowstwo(Z_\ell) \supsetneq S_0 \cap \rowstwo(Z_\ell)$, contradicting the definition of $S_0$.

Let $i \in S_0 \cap \rowstwo(Z_\ell)$. Since $i \in \rowstwo(Z_\ell)$, the function $g_i$ depends on more than one coordinate. By construction, $\row(W_0,i) \in \imsupp(\hat{g}_i)$. Since $g_i$ is Boolean-like, $|\row(W_0,i)| \ge 2$. Therefore $W_0$ contains a point $(i,j_1)$ for some $j_1 \neq j_0$. This point remains in $W_1$, showing that $i \in \rows(W_1)$. Therefore $\rows(W_1) \cap \rowstwo(Z_\ell) \supseteq S_0 \cap \rowstwo(Z_\ell)$.

By construction, $(i_1,j_0) \in W_1$, and so $i_1 \in \rows(W_1)$.
Since $i_1 \notin S_0$, this shows that $\rows(W_1) \cap \rowstwo(Z_\ell) \supsetneq S_0 \cap \rowstwo(Z_\ell)$, and so we reach a contradiction.

\paragraph{Second part} Among all sets in $\supp(\hat{f}_0)$ which agree with $S$ outside of $\rows(Z_\ell)$ and contain $\rowstwo(Z_\ell)$, let $S_1$ maximize $|S_1 \cap S \cap \rowsone(Z_\ell)|$. We will show that $S_1 \cap \rowsone(Z_\ell) \supseteq S \cap \rowsone(Z_\ell)$, and so $S_1 \supseteq S$. Since $S$ is inclusion-maximal, in fact $S_1 = S$, and so $S$ contains $\rowstwo(Z_\ell)$.

Assume, for the sake of contradiction, that $S_1 \cap \rowsone(Z_\ell)$ doesn't contain $S \cap \rowsone(Z_\ell)$. Then $S \setminus S_1$ intersects $\rows(Z_\ell)$. Since $S_1$ contains $\rowstwo(Z_\ell)$ and $\rowstwo(Z_\ell)$ is non-empty, the set $S \setminus S_1$ doesn't contain all of $\rowstwo(Z_\ell)$. Therefore \reachability shows that there exist $i_0 \in S \setminus S_1$ and $i_1 \in \rowstwo(Z_\ell)$ such that $i_0,i_1$ appear on the same column $j_0 \in \cols(Z_\ell)$.

Define $W \subseteq [n] \times [m]$ as follows: $\rows(W) = S$; for $i \in S \cap \{i_0,i_1\}$, $\row(W,i)$ is a set in $\imsupp(\hat{g}_i)$ which contains $j_0$; and $\row(W,i)$ is an arbitrary set in $\imsupp(\hat{g}_i)$ for $i \in S$ other than $i_0,i_1$.
Define $W_1$ similarly, but with $\rows(W_1) = S_1$. Note $(i_0,j_0) \in W$ (since $i_0 \in S$) and $(i_1,j_0) \in W_1$ (since $i_1 \in \rowstwo(Z_\ell) \subseteq S_1$).

Let $W_2$ be formed from $W_1$ by replacing column $j_0$ with $\col(W,j_0)$. \Duality (applied to $W,W_1,W_2$) shows that $\rows(W_2) \in \supp(\hat{f}_0)$. We will show that $\rows(W_2)$ agrees with $S_1$ outside of $\rows(Z_\ell)$, contains $\rowstwo(Z_\ell)$, and $\rows(W_2) \cap S \cap \rowsone(Z_\ell) \supsetneq S_1 \cap S \cap \rowsone(Z_\ell)$. This contradicts the definition of $S_1$.

The proof of the first two statements is identical to the proof in the first part, so we skip it.

Now suppose $i \in S_1 \cap S \cap \rowsone(Z_\ell)$, say $(i,j) \in Z$. Note $\row(W,i) = \row(W_1,i) = \{j\}$, hence both $\col(W,j)$ and $\col(W_1,j)$ contain $i$. Considering the two cases $j = j_0$ and $j \neq j_0$, we see that $i \in \col(W_2,j)$, and so $i \in \rows(W_2)$. This shows that $\rows(W_2) \cap S \cap \rowsone(Z_\ell) \supseteq S_1 \cap S \cap \rowsone(Z_\ell)$.

Since $(i_0,j_0) \in W$, by construction $i_0 \in \rows(W_2)$. By construction $i_0 \notin S_1$, and so $\rows(W_2) \cap S \cap \rowsone(Z_\ell) \supsetneq S_1 \cap S \cap \rowsone(Z_\ell)$, and we reach a contradiction.
\end{proof}

\section{Inclusion-maximal support of \texorpdfstring{$f_j,g_i$}{inner functions}}
\label{sec:step2}

In this section we prove the following result.

\begin{proposition} \label{pro:step2}
Let $f_0,f_1,\dots,f_m,g_0,g_1,\dots,g_n$ be a generalized polymorphism, where $f_0,g_0$ depend on all coordinates, and $f_1,\dots,f_m,g_1,\dots,g_n$ are Boolean-like and non-constant.

Let $i \in [n]$. Then all sets in $\imsupp(\hat{g}_i)$ contain $\rowtwo(i)$.
\end{proposition}
\begin{proof}
We start by centering the functions. For all $i \in [n]$, define $G_i = g_i - \EE[g_i]$, and for all $j \in [m]$, define $F_j = f_j - \EE[f_j]$. Let
\begin{align*}
F_0(x) &= f_0(x_1 + \EE[g_1], \dots, x_n + \EE[g_n]), &
G_0(y) &= g_0(y_1 + \EE[f_1], \dots, y_m + \EE[f_m]).
\end{align*}
Then $F_0,F_1,\dots,F_n,G_0,G_1,\dots,G_m$ is a generalized polymorphism, where $F_0,G_0$ depend on all coordinates and $F_1,\dots,F_n,G_1,\dots,G_m$ are Boolean-like (by \Cref{obs:boolean-like}), centered, and non-zero. Also, $\imsupp(\hat{g}_i) = \imsupp(\hat{G}_i)$ for all $i$, by \Cref{obs:shifting}.

\smallskip

Suppose, for the sake of contradiction, that $i_0 \in [n]$ and the set $T_{i_0} \in \imsupp(\hat{G}_{i_0})$ doesn't contain $\rowtwo(i_0)$. Let $Z_\ell$ be the connected component that contains row~$i_0$.

We claim that $i_0 \in \rowstwo(Z_\ell)$.
Indeed, if $|\row(i_0)| = 1$ then either $\rowtwo(i_0) = \emptyset$, in which case $T_{i_0}$ contains $\rowtwo(i_0)$, or $\rowtwo(i_0) = \row(i_0)$, in which case $T_{i_0} = \row(i_0)$, and in both cases we reach a contradiction.
%Hence we can assume that $|\row(i_0)| \ge 2$. %Since $G_{i_0}$ is Boolean-like, this means that $|T_{i_0}| \ge 2$. % Do we need this?

Since $T_{i_0}$ doesn't contain $\rowtwo(i_0)$, there exists some $j_0 \in \rowtwo(i_0) \setminus T_{i_0}$. By definition, $|\col(j_0)| \ge 2$, and so $\col(j_0)$ contains some $i_1 \neq i_0$. Note that $i_1 \in \rows(Z_\ell)$.

Since $f_0$ depends on all coordinates, in particular it depends on $i_1$. Therefore we can find $S \in \imsupp(\hat{F}_0)$ which contains $i_1$. Since $S$ intersects $\rows(Z_\ell)$, according to \Cref{pro:step1} it contains all of $\rowstwo(Z_\ell)$, and in particular, it contains $i_0$.

Define $W \subseteq [n] \times [m]$ as follows: $\rows(W) = S$; $\row(W,i_0) = T_{i_0}$; $\row(W,i_1)$ is some set in $\imsupp(\hat{g}_{i_1})$ containing $j_0$ (this is possible since $(i_1,j_0) \in Z$); and $\row(W,i)$ is an arbitrary set in $\imsupp(\hat{g}_i)$ for $i \in S$ other than $i_0,i_1$.
Applying \duality, we deduce that $\cols(W) \in \imsupp(\hat{G}_0)$ and $\col(W,j) \in \imsupp(\hat{F}_j)$ for all $j \in \cols(W)$. Note $(i_1,j_0) \in W$.

Let $W'$ be obtained from $W$ by replacing column $j_0$ with some set in $\imsupp(\hat{f}_{j_0})$ which contains $i_0$; this is possible since $(i_0,j_0) \in Z$. Applying \duality again, we deduce that $\row(W',i_0) \in \supp(\hat{G}_i)$. Since $j_0 \notin T_{i_0} = \row(W,i_0)$ whereas $j_0 \in \row(W',i_0)$, we have $\row(W',i_0) = T_{i_0} \cup \{j_0\}$, which contradicts the inclusion-maximality of $T_{i_0}$.
\end{proof}

\section{Shifting and structure of \texorpdfstring{$f_j,g_i$}{inner functions}}
\label{sec:step3}

In this section we prove the following result.

\begin{proposition} \label{pro:step3}
Let $f_0,f_1,\ldots,f_m,g_0,g_1,\ldots,g_n$ be a Boolean generalized polymorphism, where $f_0,g_0$ depend on all coordinates, and $f_1,\dots,f_m,g_1,\dots,g_n$ are non-constant.

There are constants $\kappa_{ij} \in \RR$ for each $i \in [n]$ and $j \in [m]$, $D_i \in \RR$ for each $i \in [n]$, and $B_j \in \RR$ for each $j \in [m]$, such that the functions
\begin{align*}
    G_i(y_1,\dots,y_m) &= g_i(y_1 - \kappa_{i1}, \dots, y_m - \kappa_{im}) - D_i, &
    F_j(x_1,\dots,x_n) &= f_j(x_1 - \kappa_{1j}, \dots, x_n - \kappa_{nj} - B_j
\end{align*}
are Boolean-like, centered, and non-zero. Together with the functions
\begin{align*}
    F_0(x_1,\dots,x_n) &= f_0(x_1 + D_1, \dots, x_n + D_n), &
    G_0(y_1,\dots,y_m) &= g_0(y_1 + B_1, \dots, y_m + B_m),
\end{align*}
which depend on all coordinates, they form a (possibly non-Boolean) generalized polymorphism.

Furthermore, if $\row(i) = \{ j \}$ then $\kappa_{ij} = D_i = 0$ (implying $G_i = g_i$), and otherwise there exists a function $\gamma_i\colon \pmone^{\rowone(i)} \to \RR$ such that
\[
 G_i(y) = \gamma_i(y|_{\rowone(i)}) \prod_{j \in \rowtwo(i)} y_j.
\]

Similarly, if $\col(j) = \{ i \}$ then $\kappa_{ij} = B_j = 0$ (implying $F_j = f_j$), and otherwise there exists a function $\phi_j\colon \pmone^{\colone(j)} \to \RR$ such that
\[
 F_j(x) = \phi_j(x|_{\colone(j)}) \prod_{i \in \coltwo(j)} x_i.
\]
\end{proposition}
\begin{proof}
For each $i \in [n]$, define $W_i \subseteq [n] \times [m]$ as follows: $\rows(W_i)$ is an arbitrary set in $\imsupp(\hat{f}_0)$ containing $i$ (this is possible since $f_0$ depends on all coordinates), and for each $i' \in \rows(W_i)$, $\row(W_i,i')$ is an arbitrary set in $\imsupp(\hat{g}_{i'})$.

\paragraph{Defining $\kappa_{ij}$}
We can assume that $(i,j) \in Z$, since otherwise $\kappa_{ij}$ doesn't affect $F_j,G_i$.
If $|\row(i)| = 1$ or $|\col(j)| = 1$ then we define $\kappa_{ij} = 0$. Otherwise, we will define $\kappa_{ij}$ so that $\hat{F}_j(\col(W_i,j) \setminus \{i\}) = 0$. Note $i \in \col(W_i,j)$: since $j \in \rowtwo(i)$, \Cref{pro:step2} shows that $j \in \row(W_i,j)$, hence $i \in \col(W_i,j)$.

Observe that the only supersets of $\col(W_i,j) \setminus \{i\}$ in $\supp(\hat{f}_j)$ are $\col(W_i,j)$ and, possibly, $\col(W_i,j) \setminus \{i\})$. Indeed, \duality shows that $\col(W_i,j) \in \imsupp(\hat{f}_j)$. Conversely, suppose that $S' \supseteq \col(W_i,j) \setminus \{i\}$ is in $\supp(\hat{f}_j)$, and choose a superset $S'' \supseteq S'$ in $\imsupp(\hat{f}_j)$. Since $i \in \coltwo(j)$ by assumption, \Cref{pro:step2} shows that $i \in S''$, and so $S'' \supseteq \col(W_i,j)$. Since $\col(W_i,j)$ is inclusion-maximal, in fact $S'' = \col(W_i,j)$, and so $\col(W_i,j) \setminus \{i\} \subseteq S' \subseteq \col(W_i,j)$.

This observation shows that
\[
 \hat{F}_j(\col(W_i,j) \setminus \{i\}) =
 \hat{f}_j(\col(W_i,j) \setminus \{i\}) - \kappa_{ij} \hat{f}_j(\col(W_i,j)).
\]
We choose accordingly
\[
 \kappa_{ij} = \frac{\hat{f}_j(\col(W_i,j) \setminus \{i\})}{\hat{f}_j(\col(W_i,j))},
\]
where the denominator is non-zero since $\col(W_i,j) \in \supp(\hat{f}_j)$.
This choice results in
\begin{equation} \label{eq:zero} \tag{$\ast$}
 \hat{F}_j(\col(W_i,j) \setminus \{i\}) = 0,
\end{equation}
which holds whenever $|\row(i)|,|\col(j)| \ge 2$.

\paragraph{Defining $D_i,B_j$} We define $D_1,\dots,D_n,B_1,\dots,B_m$ in the only way that makes $g_1,\dots,g_n,f_1,\dots,f_m$ centered. If $\row(i) = \{j\}$ then since $g_i$ is Boolean, the function $g_i(y) = \pm y_j$ is already centered, and so $D_i = 0$. Similarly, if $\col(j) = \{i\}$ then $f_j(x) = \pm x_i$ is already centered, and so $B_j = 0$.

\paragraph{Structure of $G_i$} Let $i \in [n]$ be such that $|\row(i)| \ge 2$. In order to show that $G_i$ has the claimed structure, we show that every set in $\supp(\hat{G}_i)$ contains $\rowtwo(i)$. This implies the existence of a multilinear $\gamma_i$ having the required property, which we can interpret as a real-valued function on $\pmone^{\rowone(i)}$.

At this point we remind the reader of \Cref{obs:shifting}, which shows that $\rowtwo(i)$ has the same meaning for the original polymorphism $f_0,f_1,\dots,f_m,\allowbreak g_0,g_1,\dots,g_n$ and for the shifted polymorphism $F_0,F_1,\dots,F_m,\allowbreak G_0,G_1,\dots,G_n$.

Suppose, for the sake of contradiction, that $T_i \in \supp(\hat{G}_i)$, yet $T_i$ doesn't contain $\rowtwo(i)$, say $j \in \rowtwo(i) \setminus T_i$. Recall that $i \in \rows(W_i)$ by construction. Let $W'_i$ be obtained from $W_i$ by replacing row $i$ with $T_i$; in particular, $(i,j) \notin W'_i$.

We claim that $j \in \cols(W'_i)$. Indeed, \duality shows that $\col(W_i,j) \in \imsupp(\hat{F}_j)$. Since $j \in \rowtwo(i)$, we have $|\col(j)| \ge 2$, and so $|\col(W_i,j)| \ge 2$ since $\hat{F}_j$ is Boolean-like (by \Cref{obs:boolean-like}). Since $W'_i$ is obtained by changing a single row, it is still the case that $j \in \cols(W'_i)$.

\Duality shows that $\col(W'_i,j) \in \supp(\hat{F}_j)$. Yet $\col(W'_i,j) = \col(W_i,j) \setminus \{i\}$, contradicting \Cref{eq:zero}.

\paragraph{Structure of $F_j$} Let $j \in [m]$ be such that $|\col(j)| \ge 2$. In order to show that $F_j$ has the claimed structure, we show that every set in $\supp(\hat{F}_j)$ contains $\coltwo(j)$.

Suppose, for the sake of contradiction, that $S_j \in \supp(\hat{F}_j)$, yet $S_j$ doesn't contain $\coltwo(j)$, say $i \in \coltwo(j) \setminus S_j$. Observe that $j \in \cols(W_i)$: by assumption, $j \in \rowtwo(i)$, and so $j \in \row(W_i,i)$ by \Cref{pro:step2}. Let $W'_i$ be obtained from $W_i$ by replacing column $j$ with $S_j$; in particular, $(i,j) \notin W'_i$.

We claim that $i \in \rows(W'_i)$. Indeed, $|\row(i)| \ge 2$ since $i \in \coltwo(j)$, and so $|\row(W_i,i)| \ge 2$ since $G_i$ is Boolean-like (by \Cref{obs:boolean-like}). Since $W'_i$ is obtained by changing a single column, it is still the case that $i \in \rows(W'_i)$. 

\Duality shows that $\row(W'_i,i) \in \supp(\hat{G}_i)$. Since $(i,j) \notin W'_i$, we have $j \notin \row(W'_i,i)$. However, this contradicts the structure of $G_i$, since $j \in \rowtwo(i)$.
\end{proof}

\section{Structure of \texorpdfstring{$f_0,g_0$}{outer functions}}
\label{sec:step4}

In this section we prove the following result.

\begin{proposition} \label{pro:step4}
Let $f_0,f_1,\dots,f_m,g_0,g_1,\dots,g_n$ be a Boolean generalized polymorphism, where $f_0,g_0$ depend on all coordinates, and $f_1,\dots,f_m,g_1,\dots,g_n$ are non-constant.

Assume further that if any of $f_1,\dots,f_m,g_1,\dots,g_n$ depends on a single coordinate, then it is a dictator rather than an anti-dictator (that is, of the form $f_i(x) = x_j$ rather than of the form $f_i(x) = -x_j$).

The functions $F_0,G_0$ defined in \Cref{pro:step3} have the following form. There exists a multilinear polynomial $H\colon \RR^k \to \RR$ (where $k$ is the number of connected components of $Z$) such that
\begin{align*}
F_0(x) &= H\bigl(F_0^{(1)}(x|_{\rows(Z_1)}), \dots, F_0^{(k)}(x|_{\rows(Z_k)})\bigr), &    
G_0(y) &= H\bigl(G_0^{(1)}(y|_{\cols(Z_1)}), \dots, G_0^{(k)}(y|_{\cols(Z_k)})\bigr),
\end{align*}
where for every $\ell \in [k]$, the functions $F_0^{(\ell)}$ and $G_0^{(\ell)}$ have the following form.

If $Z_\ell = \{(i,j)\}$ then $F_0^{(\ell)}(x) = x_i$ and $G_0^{(\ell)}(y) = y_j$.

If $|Z_\ell| \ge 2$ then
\begin{align*}
F_0^{(\ell)}(x) &= \prod_{i \in \rowstwo(Z_\ell)} x_i \cdot \prod_{j \in \colstwo(Z_\ell)} \phi_j(x|_{\colone(j)}), &
G_0^{(\ell)}(y) &= \prod_{j \in \colstwo(Z_\ell)} y_j \cdot \prod_{i \in \rowstwo(Z_\ell)} \gamma_i(y|_{\rowone(i)}).
\end{align*}
\end{proposition}
\begin{proof}
We first observe that all sets in $\supp(\hat{F}_0)$ intersecting $\rows(Z_\ell)$ contain $\rowstwo(Z_\ell)$. Indeed, suppose that $S_0 \in \supp(\hat{F}_0)$ intersects $\rows(Z_\ell)$ but does not contain $\rowstwo(Z_\ell)$. \Reachability shows that there exist $i_0 \in S_0 \cap \rows(Z_\ell)$ and $i_1 \in \rowstwo(Z_\ell) \setminus S_0$ that appear on the same column $j \in \cols(Z_\ell)$.

Define $W \subseteq [n] \times [m]$ as follows: $\rows(W) = S_0$; $\row(W,i_0)$ is some set in $\imsupp(\hat{G}_{i_0})$ containing $j$; and $\row(W,i)$ is an arbitrary set in $\imsupp(\hat{G}_i)$ for $i \in S_0$ other than $i_0$. \Duality shows that $\col(W,j) \in \imsupp(\hat{F}_j)$. Yet by construction, $i_1 \notin \col(W,j)$, contradicting \Cref{pro:step3}.

This (together with the dual statement about $\hat{G}_0$) allows us to expand $\hat{F}_0,\hat{G}_0$ as follows, where we separate each set of components $K \subseteq [k]$ into the set $\Kone$ of singleton components and the set $\Ktwo$ of non-singleton components:
\begin{equation} \label{eq:step4-expansion}
\begin{aligned}
F_0(x) &= \sum_{K \subseteq [\ell]}
\prod_{\substack{\ell \in \Kone \\ Z_\ell = \{(i,j)\}}} x_i \cdot
\prod_{\substack{\ell \in \Ktwo \\ i \in \rowstwo(Z_\ell)}} x_i \cdot
\Phi^K(x|_{\rowsone(Z_\ell)} : \ell \in \Ktwo), \\
G_0(y) &= \sum_{K \subseteq [\ell]}
\prod_{\substack{\ell \in \Kone \\ Z_\ell = \{(i,j)\}}} y_j \cdot
\prod_{\substack{\ell \in \Ktwo \\ j \in \colstwo(Z_\ell)}} y_j \cdot
\Gamma^K(y|_{\colsone(Z_\ell)} : \ell \in \Ktwo),
\end{aligned}
\end{equation}
for some arbitrary functions $\Phi^K,\Gamma^K$ depending on the stated coordinates.

We now substitute these expressions, and the expressions for $F_1,\dots,F_m,G_1,\dots,G_n$ stated in \Cref{pro:step3}, inside the definition of a generalized polymorphism (\Cref{eq:generalized-polymorphism}). Since $F_1,\dots,F_m,G_1,\dots,G_n$ are centered, we can trace each $z$-monomial in the left-hand side of the substituted \Cref{eq:generalized-polymorphism} to the $x$-monomial in \Cref{eq:step4-expansion} that generated it, and so to the corresponding value of $K$. The same holds for the right-hand side, and so for each $K \subseteq [k]$, equating the corresponding parts of the substituted \Cref{eq:generalized-polymorphism}, we deduce
\begin{equation} \label{eq:step4-K}
\begin{aligned}
 &
 \prod_{\substack{\ell \in \Kone \\ Z_\ell = \{(i,j)\}}} z_{ij} \cdot
 \prod_{\substack{\ell \in \Ktwo \\ i \in \rowstwo(Z_\ell)}} \left(\gamma_i(Z|_{\{i\} \times \rowone(i)}) \prod_{j \in \rowtwo(i)} z_{ij}\right) \cdot
 \Phi^K(Z_\ell^\Phi : \ell \in \Ktwo)
 \\ =&
 \prod_{\substack{\ell \in \Kone \\ Z_\ell = \{(i,j)\}}} z_{ij} \cdot
 \prod_{\substack{\ell \in \Ktwo \\ j \in \colstwo(Z_\ell)}} \left(\phi_j(Z|_{\colone(j) \times \{j\}}) \prod_{i \in \coltwo(j)} z_{ij}\right) \cdot
 \Gamma^K(Z_\ell^\Gamma : \ell \in \Ktwo),
\end{aligned}
\end{equation}
where the inputs $Z_\ell^\Phi,Z_\ell^\Gamma$ are as illustrated in \Cref{fig:Phi-Gamma}. Here we used the assumption that functions $f_j,g_i$ depending on a single coordinate are dictators.

\begin{figure}
\centering
\begin{tikzpicture}[scale=0.5]
\tikzZ

\node[align=center] at (1.5,-1.5) {$\Ztwo$};

\draw[decoration={brace},decorate] (\xleft,0.2) -- (\xright,0.2) node [midway,above] {$Z_\ell^\Gamma$};
\draw[decoration={brace},decorate] (-0.2,-\ybottom) -- (-0.2,-\ytop) node [midway,left] {$Z_\ell^\Phi$};
\end{tikzpicture}

\caption{Coordinates of $Z_\ell$ that $\Phi^K,\Gamma^K$ depend on (if $\ell \in K$)}
\label{fig:Phi-Gamma}    
\end{figure}
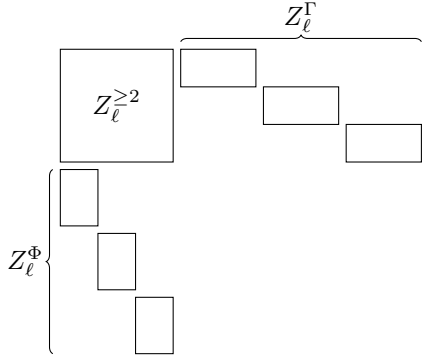

The first factor on either side of \Cref{eq:step4-K} cancels. Inside the second factor, we can cancel the product over $z_{ij}$, since on both sides, we are multiplying over all $z_{ij}$ such that $(i,j) \in \Ztwo$. Therefore for each $K \subseteq [k]$,
\begin{equation} \label{eq:step4-K-cancelled}
 \prod_{\substack{\ell \in \Ktwo \\ i \in \rowstwo(Z_\ell)}} \gamma_i(Z|_{\{i\} \times \rowone(i)}) \cdot
 \Phi^K(Z_\ell^\Phi : \ell \in \Ktwo)
 =
 \prod_{\substack{\ell \in \Ktwo \\ j \in \colstwo(Z_\ell)}} \phi_j(Z|_{\colone(j) \times \{j\}}) \cdot
 \Gamma^K(Z_\ell^\Gamma : \ell \in \Ktwo).
\end{equation}
Observe that the first factor on the left-hand side depends on $Z_\ell^\Gamma$, and similarly, the first factor on the right-hand side depends on $Z_\ell^\Phi$.

Let $\ell$ be such that $|Z_\ell| \ge 2$. For each $i \in \rowstwo(Z_\ell)$, the function $G_i$ is non-constant, and in particular, $\gamma_i$ is non-zero. Therefore we can find an assignment $U$ to $Z_\ell^\Gamma$ such that $\gamma_i(U|_{\{i\} \times \rowone(i)}) \neq 0$ for all $i \in \rowstwo(Z_\ell)$. Substituting this assignment in \Cref{eq:step4-K-cancelled} and rearranging, we obtain
\[
 \Phi^K(Z_\ell^\Phi : \ell \in \Ktwo)
 = H_{\Phi^K} \cdot \prod_{\substack{\ell \in \Ktwo \\ j \in \colstwo(Z_\ell)}} \phi_j(Z|_{\colone(j) \times \{j\}}), \text{ where }
 H_{\Phi^K} = \frac{
 \Gamma^K(U_\ell^\Gamma : \ell \in \Ktwo)
 }{
 \displaystyle\prod_{\substack{\ell \in \Ktwo \\ i \in \rowstwo(Z_\ell)}} \gamma_i(U|_{\{i\} \times \rowone(i)})
 },
\]
and so, undoing the substitution of dictatorial $G_i$'s,
\[
 \Phi^K(x|_{\rowsone(Z_\ell)} : \ell \in \Ktwo) =
 H_{\Phi^k} \cdot \prod_{\substack{\ell \in \Ktwo \\ j \in \colstwo(Z_\ell)}} \phi_j(x|_{\colone(j)}).
\]
Looking back at \Cref{eq:step4-expansion}, this shows that
\[
 F_0(x) = H_{F_0}(F_0^{(1)}(x), \dots, F_0^{(k)}(x)),
\]
where $H_{F_0}$ is the multilinear polynomial whose coefficients are $H_{\Phi^K}$.

Similarly, we can find an assignment $V$ to $Z_\ell^\Phi$, for all $\ell \in [k]$ such that $|Z_\ell| \ge 2$, such that
\[
 \Gamma^K(Z_\ell^\Gamma : \ell \in \Ktwo)
 = H_{\Gamma^K} \cdot \prod_{\substack{\ell \in \Ktwo \\ i \in \rowstwo(Z_\ell)}} \gamma_i(Z|_{\{i\} \times \rowone(i)}), \text{ where }
 H_{\Gamma^K} = \frac{
 \Phi^K(V_\ell^\Gamma : \ell \in \Ktwo)
 }{
 \displaystyle\prod_{\substack{\ell \in \Ktwo \\ j \in \colstwo(Z_\ell)}} \phi_j(V|_{\colone(j) \times \{j\}})
 }.
\]
Looking back at \Cref{eq:step4-expansion}, this shows that
\[
 G_0(y) = H_{G_0}(G_0^{(1)}(y), \dots, G_0^{(k)}(y)),
\]
where $H_{G_0}$ is the multilinear polynomial whose coefficients are $H_{\Gamma^K}$.

In order to complete the proof of the \namecref{pro:step4}, we show that $H_{\Phi^K} = H_{\Gamma^K}$ for all $K \subseteq [k]$, and so $h_{F_0} = h_{G_0}$. Indeed, observe that the two assignments $U,V$ considered above assign different variables, and so are compatible with a single assignment $W$. Substituting $W$ in \Cref{eq:step4-K-cancelled} and rearranging immediately shows that $H_{\Phi^K} = H_{\Gamma^K}$, and so completes the proof.
\end{proof}

\section{Interlude: Boolean structure lemma}
\label{sec:step5-prel}

In this section, we prove the following structure lemma, which will be used in the proof of \Cref{pro:step5}.

\begin{lemma} \label{lem:boolean-structure}
Let $a \ge 0$ and $b \ge 1$ be integers such that $a + b \ge 2$. Suppose that $f\colon \pmone^{a+b} \to \pmone$ depends on all coordinates and is given by
\[
 f(u,v) = \phi(u) \cdot \prod_{i=1}^b (v_i + \kappa_i) + B,
\]
where $\phi\colon \pmone^b \to \RR$.

Then one of the following cases holds:
\begin{itemize}
    \item Case 1: $\kappa_1 = \cdots = \kappa_b = B = 0$ and $\phi$ is $\pmone$-valued.
    \item Case 2: $\kappa_1,\dots,\kappa_b,B \in \pmone$, and $\phi$ is $\{0,C\}$-valued, where
    \[
     C = \frac{-2B}{\prod_{i=1}^b (2\kappa_i)}.
    \]
\end{itemize}
\end{lemma}
\begin{proof}
Suppose first that $a = 0$ (and so $b \ge 2$), and let $A = \phi()$. Since $f$ depends on all coordinates, necessarily $A \neq 0$. We think of $f$ as a function of $v$ alone. Choose $V_1,\dots,V_b \in \pmone$ such that $V_i \neq -\kappa_i$ for all $i \in [b]$. For $S \subseteq [b]$, let $V^{\oplus S}$ be obtained from $V$ by flipping the coordinates in $S$.

Observe that
\[
 f(V^{\oplus S}) = F \prod_{i \in S} \rho_i + B, \text{ where } F = A\prod_{i=1}^b (V_i + \kappa_i) \neq 0 \text{ and } \rho_i = \frac{-V_i + \kappa_i}{V_i + \kappa_i} \neq 1.
\]

Since $F \neq 0$ and $\rho_i \neq 1$, for every $i \in [b]$ we have $f(V^{\oplus \{i\}}) \neq f(V)$, and so $f(V^{\oplus \{i\}}) = -f(V)$. This implies that $\rho_i = -2B/F - 1$ is independent of $i$. Denote the common value by $\rho$.

Recall that $b \ge 2$. If $f(V^{\oplus \{1,2\}}) = f(V)$ then $\rho^2 = 1$, and so $\rho = -1$. For all $i \in [b]$ we have $-V_i + \kappa_i = -(V_i + \kappa_i)$, and so $\kappa_i = 0$. In view of the formula for $\rho$, in this case $B = 0$. This is Case~1.

If $f(V^{\oplus \{1,2\}}) = -f(V)$ then $\rho^2 = \rho$, and so $\rho = 0$. Therefore $\kappa_i = V_i \in \pmone$ for all $i \in [b]$. Since $f(V^{\oplus \{1,2\}}) = B$, we see that $B \in \pmone$. Finally, $-B = f(V) = A \prod_i (2\kappa_i) + B$, and so we are in Case~2, with $A = C$.

\medskip

Now suppose that $a \ge 1$. If $b \ge 2$ then we reduce to the preceding case. Since $f$ depends on all coordinates, there must be an input $U_1 \in \pmone^a$ such that $\phi(U_1) \neq 0$. Applying the case $a = 0$ to the function $f(u := U_1, v)$, we see that either $\kappa_1 = \cdots = \kappa_b = B = 0$ and $\phi(U_1) \in \pmone$, or $\kappa_1,\dots,\kappa_b,B \in \pmone$ and $\phi(U_1) = C$.

If $\kappa_1 = \cdots = \kappa_b = B = 0$ then $\phi(U) \neq 0$ for all $U \in \pmone^a$, since $f$ is Boolean. Therefore the foregoing shows that $\phi$ is $\pmone$-valued, which is Case~1.

Conversely, if $\kappa_1,\dots,\kappa_b,B \in \pmone$ then the foregoing shows that for each $U \in \pmone^a$, if $\phi(U) \neq 0$ then $\phi(U) = C$. Thus $\phi$ is $\{0,C\}$-valued, which is Case~2.

\medskip

Finally, suppose that $a \ge 1$ and $b = 1$. In this case, simplifying notation, we can write $f(u,v) = \phi(u) \cdot (v + \kappa) + B$. Since $f$ depends on all coordinates, there must be an input $U_1 \in \pmone^a$ such that $\phi(U_1) \neq 0$. Clearly $f(U_1,1) \neq f(U_1,-1)$, and so $f(U_1,1) = -f(U_1,-1)$, which implies that $\phi(U_1) \cdot \kappa = -B$.

If $\kappa = 0$ then also $B = 0$. In this case $f(u,v) = \phi(u) \cdot v$, and so clearly $\phi$ is $\pmone$-valued. This is Case~1.

If $\kappa \neq 0$ then the foregoing shows that for each $U \in \pmone^a$, if $\phi(U) \neq 0$ then $\phi(U) = -B/\kappa = -2B/(2\kappa)$.
In particular, $f(U_1,1) = -B/\kappa$.
Since $f$ depends on all coordinates, the function $\phi$ cannot be constant, and so there is an input $U_0 \in \pmone^a$ such that $\phi(U_0) = 0$. Since $f(U_0,1) = B$, we see that $B \in \pmone$. Since $\kappa = -B/f(U_1,1)$, we see that also $\kappa \in \pmone$.
This is Case~2.
\end{proof}

\section{Deducing the main theorem}
\label{sec:step5}

In this section, we prove the following result, and then explain how to deduce \Cref{thm:main} from it.

\begin{proposition} \label{pro:step5}
Let $f_0,f_1,\dots,f_m,g_0,g_1,\dots,g_n$ be a Boolean generalized polymorphism, where $f_0,g_0$ depend on all coordinates, and $f_1,\dots,f_m,g_1,\dots,g_n$ are non-constant.

Assume further that if any of $f_1,\dots,f_m,g_1,\dots,g_n$ depends on a single coordinate, then it is a dictator rather than an anti-dictator (that is, of the form $f_i(x) = x_j$ rather than of the form $f_i(x) = -x_j$).

There exists a function $h\colon \pmone^k \to \pmone$ (where $k$ is the number of connected components of $Z$) such that
\begin{align*}
 f_0(x) &= h(f_0^{(1)}(x|_{\rows(Z_1)}), \dots, f_0^{(k)}(x|_{\rows(Z_k)})), &
 g_0(y) &= h(g_0^{(1)}(y|_{\cols(Z_1)}), \dots, g_0^{(k)}(y|_{\cols(Z_k)})),
\end{align*}
for certain Boolean functions $f_0^{(\ell)},g_0^{(\ell)}$ for $\ell \in [k]$.

For each $\ell \in [k]$, the functions $f_0^{(\ell)},g_0^{(\ell)}$, the functions $f_j$ for $j \in \cols(Z_\ell)$, and the functions $g_i$ for $i \in \rows(Z_\ell)$, are as follows.

If $j \in \colsone(Z_\ell)$ and $(i,j) \in Z_\ell$ then $f_j(x) = x_i$.

If $i \in \rowsone(Z_\ell)$ and $(i,j) \in Z_\ell$ then $g_i(y) = y_j$.

For the remaining functions, we have three different possibilities.

\paragraph{Singleton case} $Z_\ell = \{(i,j)\}$, $f_0^{(\ell)}(x) = x_i$, and $g_0^{(\ell)}(y) = y_j$.

\paragraph{XOR case} There are functions $\gamma_i\colon \pmone^{\rowone(i)} \to \pmone$ for each $i \in \rowstwo(Z_\ell)$ and $\phi_j\colon \pmone^{\colone(j)} \to \pmone$ for each $j \in \colstwo(Z_\ell)$ such that
\begin{align*}
 f_0^{(\ell)}(x) &= \prod_{i \in \rowstwo(Z_\ell)} x_i \cdot \prod_{j \in \colstwo(Z_\ell)} \phi_j(x|_{\colone(j)}), &
 g_i(y) &= \prod_{j \in \rowtwo(i)} y_j \cdot \gamma_i(y|_{\rowone(i)}), \\
 g_0^{(\ell)}(y) &= \prod_{j \in \colstwo(Z_\ell)} y_j \cdot \prod_{i \in \rowstwo(Z_\ell)} \gamma_i(y|_{\rowone(i)}), &
 f_j(x) &= \prod_{i \in \coltwo(j)} x_i \cdot \phi_j(x|_{\colone(j)}).
\end{align*}

\paragraph{AND case} There are functions $\gamma_i,\phi_j$ as in the XOR case;
%There are functions $\gamma_i\colon \{0,1\}^{\rowone(i)} \to \{0,1\}$ for each $i \in \rowstwo(Z_\ell)$ and $\phi_j\colon \{0,1\}^{\colone(j)} \to \{0,1\}$ for each $j \in \colstwo(Z_\ell)$;
constants $D_i \in \pmone$ for each $i \in \rowstwo(Z_\ell)$ and $B_j \in \pmone$ for each $j \in \colstwo(Z_\ell)$; and constants $\kappa_{ij} \in \pmone$ for all $(i,j) \in \Ztwo$, such that
\begin{align*}
f_0^{(\ell)}(x) &=
2\prod_{i \in \rowstwo(Z_\ell)} \frac{x_i - D_i}{-2D_i}
\prod_{j \in \colstwo(Z_\ell)} \frac{\phi_j(x|_{\colone(j)}) + 1}{2} - 1, \\
g_i(y) &= -2D_i \frac{\gamma_i(y|_{\rowsone(i)}) + 1}{2} \prod_{j \in \rowtwo(i)} \frac{y_j + \kappa_{ij}}{2\kappa_{ij}} + D_i, \\
g_0^{(\ell)}(y) &=
2\prod_{j \in \colstwo(Z_\ell)} \frac{y_j - B_j}{-2B_j}
\prod_{i \in \rowstwo(Z_\ell)} \frac{\gamma_i(y|_{\rowone(i)}) + 1}{2} - 1, \\
f_j(x) &= -2B_j \frac{\phi_j(x|_{\colsone(i)}) + 1}{2} \prod_{i \in \coltwo(j)} \frac{x_i + \kappa_{ij}}{2\kappa_{ij}} + B_j.
\end{align*}
\end{proposition}
\begin{proof}
Let $F_0,F_1,\dots,F_m,G_1,\dots,G_n$ be the shifted generalized polymorphism constructed in \Cref{pro:step3}. Expressing the original generalized polymorphism in terms of the shifted one and using the structure given by \Cref{pro:step3,pro:step4}, we obtain that the functions $f_0,f_1,\dots,f_m,g_0,g_1,\dots,g_n$ have the structure described below.

We start with the following observation: $\kappa_{ij} = 0$
unless $|\row(i)|,|\col(j)| \ge 2$. Similarly, $D_i = 0$ unless $|\row(i)| \ge 2$, and $B_j = 0$ unless $|\col(j)| \ge 2$.

If $j \in \colsone(Z_\ell)$ and $(i,j) \in Z_\ell$ then $f_j(x) = x_i$. Otherwise,
\[
 f_j(x) = \phi_j(x|_{\colone(j)}) \prod_{i \in \coltwo(j)} (x_i + \kappa_{ij}) + B_j.
\]
Similarly, if $i \in \rowsone(Z_\ell)$ and $(i,j) \in Z_\ell$ then $g_i(y) = y_j$. Otherwise,
\[
 g_i(y) = \gamma_i(y|_{\rowone(i)}) \prod_{j \in \rowtwo(i)} (y_j + \kappa_{ij}) + D_i.
\]

Moreover, there exists a multilinear polynomial $h\colon \RR^k \to \RR$ such that
\begin{align*}
f_0(x) &= h(f_0^{(1)}(x|_{\rows(Z_1)}), \dots, f_0^{(k)}(x|_{\rows(Z_k)})), &    
f_0(y) &= h(g_0^{(1)}(y|_{\cols(Z_1)}), \dots, g_0^{(k)}(y|_{\cols(Z_k)})),
\end{align*}
where the functions $f_0^{(\ell)},g_0^{(\ell)}$ have the following form. If $Z_\ell = \{(i,j)\}$ then $f_0^{(\ell)}(x) = x_i$ and $g_0^{(\ell)}(y) = y_j$, and if $|Z_\ell| \ge 2$ then
\begin{align*}
f_0^{(\ell)}(x) &= \prod_{i \in \rowstwo(Z_\ell)} (x_i - D_i) \cdot \prod_{j \in \colstwo(Z_\ell)} \phi_j(x|_{\colone(j)}), &
g_0^{(\ell)}(y) &= \prod_{j \in \colstwo(Z_\ell)} (y_j - B_j) \cdot \prod_{i \in \rowstwo(Z_\ell)} \gamma_i(y|_{\rowone(i)}).
\end{align*}

\medskip

The above already handles connected components which are singletons.
In order to handle the remaining components, we need (in some cases) to apply (invertible) affine shifts to the functions $\gamma_i,\phi_j,f_0^{(\ell)},g_0^{(\ell)}$, and then to update $h$ accordingly, so that the formulas for $f_0,g_0$ still hold; this will be possible if the affine shift applied to $f_0^{(\ell)}$ and $g_0^{(\ell)}$ is the same. We denote the new versions of these functions by tildes. If no shift is needed for a connected components (which is the case for singleton connected components, for example), we apply the identity shift.

\smallskip

Next, we handle connected components whose points are all on the same row or column. Suppose that $Z_\ell$ is a non-singleton components all of whose points are on row $i$ (the case in which all points are on the same column is analogous). Then $f_j(x) = x_i$ for all $j \in \cols(Z_\ell)$, and
\begin{align*}
g_i(y) &= \gamma_i(y|_{\row(i)}) + D_i, &
f_0^{(\ell)}(x) &= x_i - D_i, &
g_0^{(\ell)}(y) &= \gamma_i(y|_{\row(i)}).
\end{align*}
We shift these functions as follows:
\begin{align*}
\tilde \gamma_i &= \gamma_i + D_i, &
\tilde f_0^{(\ell)} &= f_0^{(\ell)} + D_i = x_i, &
\tilde g_0^{(\ell)} &= g_0^{(\ell)} + D_i = \tilde \gamma_i,
\end{align*}
so that $g_i = \tilde\gamma_i$. The functions $\tilde\gamma_i,\tilde f_0^{(\ell)},\tilde g_0^{(\ell)}$ are Boolean, and so this conforms with the XOR case.

\smallskip

Finally, suppose that $Z_\ell$ is a non-singleton component in which not all points are on the same row or column. This implies that for each $i \in \rows(Z_\ell)$ we have $\rowtwo(i) \neq \emptyset$, and similarly for each $j \in \cols(Z_\ell)$ we have $\coltwo(i) \neq \emptyset$.
At this point, we bring into the fold \Cref{lem:boolean-structure}, noting that the lemma applies to the expressions for $g_i$ (whenever $|\row(i)| \ge 2$) and $f_j$ (whenever $|\col(j)| \ge 2$). Indeed, considering $g_i$, we have $a = |\rowone(i)|$ and $b = |\rowtwo(i)|$. As seen above, $b \ge 1$, and by assumption, $a + b = |\row(i)| \ge 2$.

For every $i \in \rowstwo(Z_\ell)$ and every $j \in \colstwo(Z_\ell)$, \Cref{lem:boolean-structure} gives us one of two possible structures, and these can be told apart by considering $\kappa_{ij}$ for any $j \in \rowtwo(i)$ (in the case of $g_i$) or any $i \in \coltwo(j)$ (in the case of $f_j$). Since $\Ztwo$ is itself connected (as a subgraph of the graph we defined on $Z$ back in \Cref{sec:introduction}), the same structure must hold for \emph{all} $i \in \rowstwo(Z_\ell)$ and $j \in \colstwo(Z_\ell)$.

We now split into two cases, depending on whether Case~1 holds or Case~2 holds. If Case~1 holds then $\kappa_{ij} = 0$ for all $(i,j) \in \Ztwo$, $D_i = 0$ for all $i \in \rowstwo(Z_\ell)$, and $B_j = 0$ for all $j \in \colstwo(Z_\ell)$. Moreover, the functions $\gamma_i,\phi_j$ are Boolean. This conforms to the XOR case.

If Case~2 holds then $\kappa_{ij} \in \pmone$ for all $(i,j) \in \Ztwo$, $D_i \in \pmone$ for all $i \in \rowstwo(Z_\ell)$, and $B_j \in \pmone$ for all $j \in \colstwo(Z_\ell)$. Moreover, for each $i \in \rowstwo(Z_\ell)$, the function $\gamma_i$ is $\{0,C_i\}$-valued, where $C_i = -2D_i/\prod_{j \in \rowtwo(i)} (2\kappa_{ij})$. Similarly, for each $j \in \colstwo(Z_\ell)$, the function $\phi_j$ is $\{0,A_j\}$-valued, where $A_j = -2B_j/\prod_{i \in \coltwo(j)} (2\kappa_{ij})$.

In order to make the functions $\gamma_i,\phi_j$ Boolean, we apply the following affine shifts:
\begin{align*}
\tilde\gamma_i &= 2\gamma_i/C_i - 1, &
\tilde\phi_j &= 2\phi_j/A_j - 1,
\end{align*}
and so
\begin{align*}
g_i(y) &= \frac{C_i}{2} (\tilde\gamma_i(y|_{\rowsone(i)}) + 1) \prod_{j \in \rowtwo(i)} (y_j + \kappa_{ij}) + D_i =
-2D_i \frac{\tilde\gamma_i(y|_{\rowsone(i)}) + 1}{2} \prod_{j \in \rowtwo(i)} \frac{y_j + \kappa_{ij}}{2\kappa_{ij}} + D_i, \\
f_j(x) &= \frac{A_j}{2} (\tilde\phi_i(x|_{\colsone(j)}) + 1) \prod_{i \in \coltwo(j)} (x_i + \kappa_{ij}) + B_j =
-2B_j \frac{\tilde\phi_j(x|_{\colsone(i)}) + 1}{2} \prod_{i \in \coltwo(j)} \frac{x_i + \kappa_{ij}}{2\kappa_{ij}} + B_j.
\end{align*}

We shift $f_0^{(\ell)}$ and $g_0^{(\ell)}$ accordingly:
\begin{align*}
\tilde f_0^{(\ell)}(x) &= 2 \frac{\prod_{j \in \colstwo(Z_\ell)} (1/A_j)}{\prod_{i \in \rowstwo(Z_\ell)} (-2D_i)} f_0^{(\ell)} - 1 =
2\prod_{i \in \rowstwo(Z_\ell)} \frac{x_i - D_i}{-2D_i}
\prod_{j \in \colstwo(Z_\ell)} \frac{\tilde\phi_j(x|_{\colone(j)}) + 1}{2} - 1, \\
\tilde g_0^{(\ell)}(y) &= 2 \frac{\prod_{i \in \rowstwo(Z_\ell)} (1/C_i)}{\prod_{j \in \colstwo(Z_\ell)} (-2B_j)} g_0^{(\ell)} - 1 =
2\prod_{j \in \colstwo(Z_\ell)} \frac{y_j - B_j}{-2B_j}
\prod_{i \in \rowstwo(Z_\ell)} \frac{\tilde\gamma_i(y|_{\rowone(i)}) + 1}{2} - 1.
\end{align*}
Despite the different formulas, this is the same affine shift, since
\[
\frac{\prod_{j \in \colstwo(Z_\ell)} (1/A_j)}{\prod_{i \in \rowstwo(Z_\ell)} (-2D_i)}
=
\frac{\prod_{(i,j) \in \Ztwo} (2\kappa_{ij})}{\prod_{i \in \rowstwo(Z_\ell)} (-2D_i) \prod_{j \in \colstwo(Z_\ell)} (-2B_j)}
=
\frac{\prod_{i \in \rowstwo(Z_\ell)} (1/C_i)}{\prod_{j \in \colstwo(Z_\ell)} (-2B_j)}
.
\]
By construction, the functions $\tilde f_0^{(\ell)}$ and $\tilde g_0^{(\ell)}$ are Boolean, and so this conforms to Case~2.

\medskip

Finally, since the functions $\tilde f_0^{(1)},\dots,\tilde f_0^{(k)}$ are Boolean and non-constant, we see that $\tilde h$ (considered as a function on $\pmone^k$) is Boolean.
\end{proof}

It remains to explain how to obtain the statement of \Cref{thm:main} from that of \Cref{pro:step5}. We do this in two steps. First, we interpret the formulas of \Cref{pro:step5} in terms of Boolean logic. Second, we allow anti-dictators.

In order to interpret the formulas \Cref{pro:step5} in terms of Boolean logic, we first explain how we convert $\pmone$ to $\{0,1\}$: using $(-1,1) \mapsto (1,0)$. This ensures that product corresponds to XOR. Taking $\sigma_j = 0$ and $\tau_i = 0$ for all relevant $i,j$, this explains the formulas for everything other than the AND case.

As for the AND case, let us consider $g_i$ first (the function $f_j$ is analogous). If $y_j = \kappa_{ij}$ for all $j \in \rowtwo(i)$ and $\gamma_i(y|_{\rowone(i)}) = 1$, then $g_i(y) = -D_i$, and otherwise $g_i(y) = D_i$. This becomes the AND case of \Cref{thm:main} once we flip $\kappa_{ij},\gamma_i$.

Next, consider $g_0^{(\ell)}$ (the function $f_0^{(\ell)}$ is analogous). If $y_j = -B_j$ for all $j \in \colstwo(Z_\ell)$ and $\gamma_i(y|_{\rowone(i)}) = 1$ for all $i \in \rowstwo(Z_\ell)$ then $g_0^{(\ell)}(y) = 1$, and otherwise $g_0^{(\ell)}(y) = -1$. This becomes the AND case of \Cref{thm:main} since we flipped $\gamma_i$ and once we flip $g_0^{(\ell)}$ itself (and update $h$ accordingly).

\medskip

The other aspect in which \Cref{pro:step5} and \Cref{thm:main} differ is that the latter allows anti-dictators. That is, if $\row(i) = \{j\}$ then $g_i(y) = y_j \oplus \tau_i$ for some $\tau_i \in \{0,1\}$, and if $\col(j) = \{i\}$ then $f_j(x) = x_i \oplus \sigma_j$ for some $\sigma_j \in \{0,1\}$. 

We can obtain a Boolean generalized polymorphism without anti-dictators as follows: for each row $i$ such that $|\row(i)| = 1$, replace $g_i$ by $g_i \oplus \tau_i$, and replace the $i$'th input to $f_0$ by $x_i \oplus \tau_i$; and operate analogously on columns $j$ such that $|\col(j)| = 1$. \Cref{pro:step5} (in its logical form) now applies. Undoing this change, we obtain \Cref{thm:main} as stated.

\bibliographystyle{alpha}
\bibliography{biblio}

\end{document}